\documentclass[10pt, leqno]{amsart}

\usepackage{amssymb,amsfonts,amscd,amsbsy}
\usepackage[mathscr]{eucal}

\theoremstyle{plain}
\newtheorem{thm}{Theorem}[section]

\theoremstyle{definition}

\newtheorem{rem}[thm]{Remark}
\newtheorem{exa}[thm]{Example}

\begin{document}

\title[Congruences]{Congruences for Traces of Singular moduli}
\author[Robert Osburn]{Robert Osburn}

\address{Department of Mathematics $\&$ Statistics, Queen's University, Kingston, Ontario, Canada K7L 3N6}

\email{osburnr@mast.queensu.ca}


\begin{abstract}
We extend a result of Ahlgren and Ono \cite{ao} on congruences for
traces of singular moduli of level $1$ to traces defined in terms
of Hauptmodul associated to certain groups of genus 0 of higher
levels.
\end{abstract}

\maketitle
\section{Introduction}

Let $j(z)$ denote the usual elliptic modular function on
$\operatorname{SL}_2(\mathbb Z)$ with $q$-expansion ($q:=e^{2\pi i
z}$)

\begin{center}
$j(z)=q^{-1}+744+196884q+21493760q^2+\cdots.$
\end{center}

The values of $j(z)$ at imaginary quadratic arguments in the upper
half of the complex plane are known as singular moduli. Singular
moduli are important algebraic integers which generate ring class
field extensions of imaginary quadratic fields (Theorem 11.1 in
\cite{cox}), are related to supersingular elliptic curves
(\cite{ao}), and to Borcherds products of modular forms
(\cite{b1}, \cite{b2}).

Let $d$ denote a positive integer congruent to 0 or 3 modulo 4 so
that $-d$ is the discriminant of an order in an imaginary
quadratic field. Denote by $\mathcal{Q}_d$ the set of positive
definite integral binary quadratic forms
$$
Q(x,y)=ax^2+bxy+cy^2
$$
with discriminant $-d=b^2-4ac$. To each $Q \in \mathcal{Q}_d$, let
$\alpha_Q$ be the unique complex number in the upper half plane
which is a root of $Q(x,1)$;  the singular modulus $j(\alpha_Q)$
depends only  on the equivalence class of $Q$ under the action of
$\Gamma:= PSL_2(\mathbb Z)$. Define $\omega_Q\in \{1, 2, 3\}$ by
$$
\omega_Q:
=\left \{ \begin{array}{l}
2 \quad \mbox{if $Q \sim_{\Gamma} [a,0,a]$},\\
3 \quad  \mbox{if $Q \sim_{\Gamma} [a,a,a]$},\\
1 \quad \mbox{otherwise.}
\end{array}
\right. \\
$$
Let $J(z)$ be the Hauptmodul
$$
  J(z):= j(z)-744=q^{-1}+196884q+21493760q^2+\cdots.
$$
Zagier \cite{z} defined the trace of the singular moduli of
discriminant $-d$ as
$$
   t(d):=\sum_{Q\in \mathcal{Q}_d/\Gamma}
   \frac{J(\alpha_Q)}{\omega_Q}=\sum_{Q\in \mathcal{Q}_d/\Gamma}
   \frac{j(\alpha_Q)-744}{\omega_Q} \in \mathbb{Z}.
$$

Zagier has shown that $t(d)$ has some interesting properties.
Namely, the following result (see Theorem 1 in \cite{z}) shows
that the $t(d)$'s are Fourier coefficients of a half-integral
weight modular form.

\begin{thm} Let $\theta_1(z)$ and $E_4(z)$ be defined by
$$
\begin{aligned}
&E_4(z):=1+240\sum_{n=1}^{\infty}\frac{n^3q^n}{1-q^n},\\
&\theta_1(z):=
\frac{\eta^2(z)}{\eta(2z)}
=\sum_{n=-\infty}^{\infty}(-1)^nq^{n^2}=1-2q+2q^4-2q^9+\cdots.\\
\end{aligned}
$$
and let $g(z)$ be defined by
$$
\begin{aligned}
g(z): &=-q^{-1}+2+\sum_{0<
d\equiv 0, 3\pmod 4}t(d)q^d \\
\end{aligned}
$$

Then

$$
\begin{aligned}
g(z) &=-\frac{\theta_1(z)E_4(4z)}{\eta^6(4z)} \\
&=-q^{-1}+2-248q^3+492q^4-4119q^7\cdots \\
\end{aligned}
$$

\noindent i.e., $g(z)$ is a modular form of weight $\frac{3}{2}$
on $\Gamma_{0}(4)$, holomorphic on the upper half plane and
meromorphic at the cusps.

\end{thm}

Now what about divisibility properties of $t(d)$ as $d$ varies? In
this direction, Ahlgren and Ono \cite{ao} recently proved the
following result which shows that these traces $t(d)$ satisfy
congruences based on the factorization of primes in certain
imaginary quadratic fields.

\begin{thm} If $d$ is a positive integer for which an odd prime $l$ splits in
$\mathbb Q(\sqrt{-d})$, then
$$
   t(l^2d)\equiv 0\pmod l.
$$
\end{thm}

Recently, Kim \cite{kim2} and Zagier \cite{z} defined an analogous
trace of singular moduli by replacing the $j$-function by a
modular function of higher level, in particular by the Hauptmodul
associated to other groups of genus 0. Let $\Gamma_0(N)^{*}$ be
the group generated by $\Gamma_0(N)$ and all Atkin-Lehner
involutions $W_e$ for $e||N$, i.e., $e$ is a positive divisor of $N$ for which
gcd$(e, N/e)=1$. There are only finitely many values
of $N$ for which $\Gamma_0(N)^{*}$ is of genus 0 (see
\cite{fricke1}, \cite{fricke2}, or \cite{ogg}). In particular,
there are only finitely many prime values of $N$. For such a prime
$p$, let $j_{p}^{*}$ be the corresponding Hauptmodul. For these
primes $p$, Kim and Zagier define a trace $t^{(p)}(d)$ (see
Section 3 below) in terms of singular values of $j_{p}^{*}$. The
goal of this paper is to prove that the same type of congruence
holds for $t^{(p)}(d)$, namely

\begin{thm} Let $p$ be a prime for which $\Gamma_0(p)^{*}$ is of genus 0. If $d$ is a
positive integer such that $-d$ is congruent to a square modulo
$4p$ and for which an odd prime $l \neq p$ splits in $\mathbb
Q(\sqrt{-d})$, then
$$
   t^{(p)}(l^2d)\equiv 0\bmod l.
$$
\end{thm}

\section{Preliminaries on Modular and Jacobi forms}

We first recall some facts about half-integral weight modular forms (see
\cite{kob}, \cite{koh}). If $f(z)$ is a function of the upper half-plane,
$\lambda \in \frac{1}{2}\mathbb Z$, and $\left(\begin{matrix} a & b \\
c & d \\ \end{matrix} \right) \in GL_{2}^{+}(\mathbb R)$, then we define the
slash operator by

\begin{center}
$\displaystyle f(z)|_{\lambda}\left(\begin{matrix} a & b \\
c & d \\ \end{matrix} \right) :=
(ad-bc)^{\frac{\lambda}{2}}(cz+d)^{-\lambda}f\Big(\frac{az+b}{cz+d}\Big)$
\end{center}

Here we take the branch of the square root having non-negative real part.
If $\gamma=\left(\begin{matrix} a & b \\
c & d \\ \end{matrix} \right) \in \Gamma_0(4)$, then define

\begin{center}
$\displaystyle j(\gamma, z):=\Big(\frac{c}{d}\Big)\epsilon_d^{-1}\sqrt{cz+d}$,
\end{center}

\noindent where 

$$\epsilon:=
\left \{ \begin{array}{l}
1 \quad \mbox{if $d \equiv 1 \pmod 4$},\\
i \quad  \mbox{if $d \equiv -1 \pmod 4$}.
\end{array}
\right. \\
$$

If $k$ is an integer and $N$ is an odd positive integer, then let
$\mathcal{M}_{k+{\frac{1}{2}}}(\Gamma_0(4N))$ denote the infinite
dimensional vector space of nearly holomorphic modular forms of
weight $k+\frac{1}{2}$ on $\Gamma_0(4N)$. These are functions $f(z)$ which are
holomorphic on the upper half-plane, meromorphic at the cusps, and which satisfy

\begin{equation}
f({\gamma}z)=j(\gamma, z)^{2k+1}f(z)
\end{equation}

\noindent for all $\gamma \in \Gamma_0(4N)$. Denote by
$\mathcal{M}^{+}_{k+{\frac{1}{2}}}(\Gamma_0(4N))$ the ``Kohnen
plus-spaces'' (see \cite{koh}) of nearly holomorphic forms which
transform according to (1) and which have a Fourier expansion of the form

\begin{center}
$\displaystyle \sum_{(-1)^{k}n \equiv 0,1 \pmod 4} a(n)q^n$.
\end{center}

We recall some properties of Hecke operators on
$\mathcal{M}^{+}_{k+{\frac{1}{2}}}(\Gamma_0(4N))$. If $l$ is a
prime such that $l \nmid N$, then the Hecke operator
$T_{k+\frac{1}{2}, 4N}(l^2)$ on a modular form

\begin{center}
$f(z):= \displaystyle \sum_{(-1)^{k}n \equiv 0,1 \pmod 4} a(n)q^n
\in \mathcal{M}^{+}_{k+{\frac{1}{2}}}(\Gamma_0(4N))$
\end{center}

\noindent is given by

\begin{center}
$f(z)|T_{k+\frac{1}{2}, 4N}(l^2):= \displaystyle \sum_{(-1)^{k}n
\equiv 0,1 \pmod 4} \Bigl ( a(l^{2}n) +
\Biggl(\frac{(-1)^{k}n}{l}\Biggr)l^{k-1}a(n) + l^{2k-1}a(n/l^2)
\Bigr ) q^n$
\end{center}

\noindent where $\Bigl(\frac{*}{l}\Bigr)$ is a Legendre symbol.
Let us now recall some facts about Jacobi forms (see \cite{ez}). A
Jacobi form on $\operatorname{SL}_2(\mathbb Z)$ is a holomorphic
function

\begin{center}
$\phi: \frak{H} \times \mathbb C \to \mathbb C$
\end{center}

\noindent satisfying

\begin{center}
$\displaystyle \phi\Big(\frac{a{\tau}+b}{c{\tau}+d}, \frac{z}{c{\tau} + d}\Big)=
(c{\tau}+d)^{k} e^{2{\pi}iN\frac{cz^2}{c{\tau}+ d}} \phi(\tau, z)$
\end{center}

\noindent

\begin{center}
$\displaystyle \phi({\tau}, z+{\lambda}{\tau}+ {\mu})=
e^{-2{\pi}iN({\lambda}^2{\tau}+2{\lambda}z)} \phi(\tau, z)$
\end{center}

\noindent for all $\left( \begin{matrix} a & b \\
c & d \\ \end{matrix} \right) \in \operatorname{SL}_2(\mathbb Z)$
and $(\lambda,\mu) \in {\mathbb Z}^2$, and having a Fourier
expansion of the form ($q=e^{2{\pi}i{\tau}}$,
$\zeta=e^{2{\pi}iz}$)

\begin{center}
$\phi(\tau, z) = \displaystyle \sum_{n=0}^{\infty}
\sum_{\substack{r\in \mathbb Z \\ r^2 \leq 4Nn}}
c(n,r)q^{n}{\zeta}^{r}$.
\end{center}

\noindent Here $k$ and $N$ are the weight and index of $\phi$,
respectively. Let $J_{k,N}$ denote the space of Jacobi forms of
weight $k$ and index $N$ on $\operatorname{SL}_2(\mathbb Z)$. By
Theorem 2.2 in \cite{ez}, the coefficient $c(n,r)$ depends only on
$4Nn - r^2$ and $r \bmod 2N$. By definition $c(n,r)=0$ unless $4Nn
- r^2 \geq 0$. If we drop the condition $4Nn - r^2 \geq 0$, we
obtain a nearly holomorphic Jacobi form. Let $J^{!}_{k,N}$ be the
space of nearly holomorphic Jacobi forms of weight $k$ and index
$N$.

\section{Traces}

Let $\Gamma_0(N)^{*}$ be the group generated by $\Gamma_0(N)$ and
all Atkin-Lehner involutions $W_e$ for $e||N$, that is, $e$ is a
positive divisor of $N$ for which gcd$(e, N/e)=1$. $W_e$ can be
represented by a matrix of the form $\frac{1}{\sqrt{e}}\left( \begin{matrix} ex & y \\
Nz & ew \\ \end{matrix} \right)$ with $x$, $y$, $z$, $w \in
\mathbb Z$ and $xwe-yzN/e=1$. There are only finitely many values
of $N$ for which $\Gamma_0(N)^{*}$ is of genus 0 (see \cite{fricke1},
\cite{fricke2}, or \cite{ogg}). In particular,
if we let $\frak{S}$ denote the set of prime values for such $N$,
then

\begin{center}
$\frak{S}=\{2, 3, 5, 7, 11, 13, 17, 19, 23, 29, 31, 41, 47, 59, 71
\}.$
\end{center}

For $p \in \mathfrak{S}$, let $j_{p}^{*}$ be the corresponding
Hauptmodul with Fourier expansion

\begin{center}
$q^{-1} + 0 + a_{1}q + a_{2}q^2 + \dots$
\end{center}

Let us now define a trace in terms of the $j_{p}^{*}$'s. Let $d$
be a positive integer such that $-d$ is congruent to a square
modulo $4p$. Choose an integer $\beta \bmod 2p$ such that
${\beta}^2 \equiv -d \bmod 4p$ and consider the set

\begin{center}
$\mathcal{Q}_{d, p, \beta}= \{[a,b,c] \in \mathcal{Q}_d : a \equiv
0 \bmod p$, $b \equiv \beta \bmod 2p \}$.
\end{center}

Note that $\Gamma_0(p)$ acts on $\mathcal{Q}_{d, p, \beta}$.
Assume that $d$ is not divisible as a discriminant by the square
of any prime dividing $p$, i.e. not divisible by $p^2$. Then we
have a bijection via the natural map between

\begin{center}
$\mathcal{Q}_{d, p, \beta} \diagup \Gamma_0(p)$
\end{center}

\noindent and

\begin{center}
$\mathcal{Q}_d \diagup \Gamma$
\end{center}

\noindent as the image of the root $\alpha_{Q}$, $Q \in
\mathcal{Q}_{d, p, \beta}$, in $\Gamma_0(p) \diagup \mathfrak{H}$
corresponds to a Heegner point. We could then define a trace
$t^{(p, \beta)}(d)$ as the sum of the values of $j_{p}^{*}$ with
$Q$ running over a set of representatives for $\mathcal{Q}_{d, p,
\beta} \diagup \Gamma_0(p)$. As $t^{(p, \beta)}(d)$ is independent
of $\beta$, we define the trace $t^{(p)}(d)$ (see Section 8 of
\cite{z} or Section 1 of \cite{kim2})

\begin{center}
$t^{(p)}(d) = \displaystyle \sum_{Q}
\frac{j_{p}^{*}(\alpha_{Q})}{\omega_{Q}} \in \mathbb{Z}$
\end{center}

\noindent where the sum is over $\Gamma_0(p)^{*}$ representatives
of forms $Q=[a,b,c]$ satisfying $a \equiv 0 \bmod p$.

\begin{rem}
For $p=2$, we have $t^{(2)}(4)=\frac{1}{2}
j_{2}^{*}(\frac{1+i}{2})=-52$, $t^{(2)}(7)=
j_{2}^{*}(\frac{1+\sqrt{-7}}{4})=-23$, $t^{(2)}(8)=
j_{2}^{*}(\frac{\sqrt{-2}}{2})=152$. For $p=3$, we have
$t^{(3)}(3)=\frac{1}{3} j_{3}^{*}(\frac{-3+\sqrt{-3}}{6})=-14$,
$t^{(3)}(11)= j_{3}^{*}(\frac{1+\sqrt{-11}}{6})=22$. Moreover by
the table in Section 8 of \cite{z}, we have:

\begin{center}
\begin{tabular}{c|ccc}
\multicolumn{4}{c}{} \\
$d$ & $t^{(2)}(d)$ & $t^{(3)}(d)$ & $t^{(5)}(d)$ \\
\hline 3 &  & $-14$ &  \\
\hline 4 & $-52$ &  & $-8$ \\
\hline 7 & $-23$ &  &  \\
\hline 8 & 152 & $-34$ &  \\
\hline 11 &  & 22 & $-12$ \\
\hline 12 & $-496$ & 52 &  \\
\hline 15 & $-1$ & $-138$ & $-38$ \\
\hline 16 & 1036 &  & $-6$ \\
\hline 19 &  &  & 20 \\
\hline 20 & $-2256$ & $-116$ & 12 \\
\hline 23 & $-94$ & 115 &  \\
\hline 24 & 4400 & 348 & $-44$ \\
\hline 27 &  & $-482$ &  \\
\hline 28 & $-8192$ & & \\
\end{tabular}
\end{center}

\noindent The empty entries correspond to $-d$ which are not
congruent to squares modulo $4p$.

\end{rem}

By the discussion in Section 8 of \cite{z} or Section 2.2 in
\cite{kim2}, there exist forms $\phi_{p} \in J^{!}_{2,p}$ uniquely
characterized by the condition that their Fourier coefficients
$c(n,r)=B(4pn-r^2)$ depend only on $r^2-4pn$ and where $B(-1)=1$,
$B(d)=0$ if $d=4pn-r^2<0$, $\neq -1$ and $B(0)=-2$. Define
$g_p(z)$ as

\begin{center}
$g_p(z):= q^{-1} + \displaystyle \sum_{d \geq 0} B(d)q^d$.
\end{center}

By the correspondence between Jacobi forms and half-integral
weight forms (Theorem 5.6 in \cite{ez}), $g_p(z) \in
\mathcal{M}^{+}_{\frac{3}{2}}(\Gamma_0(4p))$. As the dimension of
$J_{2,p}$ is zero, we have that for every integer $d \geq 0$ such
that $-d$ is congruent to a square modulo $4p$, there exists a
unique $f_{d,p} \in \mathcal{M}^{+}_{\frac{1}{2}}(\Gamma_0(4p))$
with Fourier expansion

\begin{center}
$f_{d,p}(z)= q^{-d} + \displaystyle \sum_{0<D\equiv 0,1 \pmod 4}
A(D,d) q^D$.
\end{center}

\noindent An explicit construction of $f_{d,p}$ can be found in
the appendix of \cite{kim1} and the uniqueness of $f_{d,p}$
follows from the discussion at the end of Section 2 in
\cite{kim1}. The following result relates the Fourier coefficients
$A(1,d)$ and $B(d)$ and shows that the traces $t^{(p)}(d)$ are
Fourier coefficients of a nearly holomorphic Jacobi form of weight
2 and index $p$ (see Theorem 8 in \cite{z} or Lemma 3.5 and
Corollary 3.6 in \cite{kim2}).

\begin{thm} Let $p$ be a prime for which $\Gamma_0(p)^{*}$ is of genus 0. \\
(i) Let $d=4pn-r^2$ for some integers $n$ and $r$. Let $A(1,d)$ be
the coefficient of $q$ in $f_{d,p}$ and $B(d)$ be the coefficient
of ${q^n}{\zeta^{r}}$ in
$\phi_{p}$. Then $A(1,d)=-B(d)$. \\
(ii) For each natural number $d$ which is congruent to a square
modulo $4p$, let $t^{(p)}(d)$ be defined as above. We also put
$t^{(p)}(-1)=-1, t^{(p)}(d)=0$ for $d<-1$. Then
$t^{(p)}(d)=-B(d)$.
\end{thm}

\section{Proof of Theorem 1.3}

\begin{proof}

The proof requires the study of Hecke operators $T_{k+\frac{1}{2},
4p}(l^2)$ on the forms $g_p(z)$ and $f_{d,p}(z)$. Define integers
$A_{l}(d)$ and $B_{l}(d)$ by

\begin{center}
$A_{l}(d):=$ the coefficient of $q$ in $f_{d,p}|T_{\frac{1}{2},
4p}(l^2)$,
\end{center}

\begin{center}
$B_{l}(d):=$ the coefficient of $q^{d}$ in $g_p(z)|T_{\frac{3}{2},
4p}(l^2)$.
\end{center}

From equation (19) of \cite{z}, we have

\begin{center}
$A_{l}(d)=A(1,d) + lA(l^2,d)$.
\end{center}

Also note that we have

\begin{center}
$g_p(z)|T_{\frac{3}{2}, 4p}(l^2)= q^{-1} + lq^{-{l^2}} +
\displaystyle \sum_{0<d \equiv 0, 3 \pmod 4} \Bigl ( B(l^{2}d) +
\Biggl(\frac{-d}{l}\Biggr)B(d) + lB(d/l^2) \Bigr ) q^d$
\end{center}

\noindent and so $B_{l}(d) = B(l^{2}d) +
\Bigl(\frac{-d}{l}\Bigr)B(d) + lB(d/l^2)$. Now suppose $p$ is in
$\frak{S}$ and $d$ is a positive integer such that $-d$ is a
square modulo $4p$ and for which an odd prime $l\neq p$ splits in
$\mathbb Q(\sqrt{-d})$. Then $\Bigl(\frac{-d}{l}\Bigr)=1$. By
Theorem 3.2 and the above calculations, we have
$$
\begin{aligned}
t^{(p)}(l^{2}d)&= -B(l^{2}d) \\
&= -B_{l}(d) + \Bigl(\frac{-d}{l}\Bigr)B(d) + lB(d/l^{2}) \\
&\equiv -B_{l}(d) + B(d) \bmod l \\
&\equiv A_{l}(d) + B(d) \bmod l \\
&\equiv A(1,d) + lA(l^{2},d) + B(d) \bmod l \\
&\equiv -B(d) + B(d) \bmod l \\
&\equiv 0 \bmod l.\\
\end{aligned}
$$
\end{proof}

\begin{exa} We now illustrate Theorem 1.3. If $p=2$ and $l=3$,
then for every non-negative integer $s$, we have

\begin{center}
$t^{(2)}(3^2(24s+23)) \equiv 0 \bmod 3$.
\end{center}

In particular, if we want to compute $t^{(2)}(207)$, then we are
interested in $\phi_{2} \in J^{!}_{2,2}$. By Theorem 9.3 in
\cite{ez} and the discussion preceding Table 8 in \cite{z},
$J^{!}_{2,2}$ is the free polynomial algebra over
\begin{center}
$\mathbb C[E_4(\tau), E_6(\tau),
{\Delta}^{-1}]\diagup({E_4(\tau)}^{3} - {E_6(\tau)}^{2})$
\end{center}

\noindent on two generators $a$ and $b$ where
$\displaystyle \Delta=\frac{{E_4(\tau)}^{3} - {E_6(\tau)}^{2}}{1728}$. The
Fourier expansions of $a$ and $b$ begin

$$
\begin{aligned}
a&=(\zeta - 2 + {\zeta}^{-1}) + (-2{\zeta}^{2} + 8{\zeta} - 12 +
8{\zeta}^{-1} - 2{\zeta}^{-2})q + ({\zeta}^{3} -12{\zeta}^2 +
39{\zeta} - 56 + \cdots)q^2 \\ &+ (8{\zeta}^3 - 56{\zeta}^2 +
152{\zeta} - 208 + \cdots)q^3 + \cdots.
\end{aligned}
$$

$$
\begin{aligned}
b&=(\zeta + 10 + {\zeta}^-1) + (10{\zeta}^2 -64{\zeta} + 108 -
64{\zeta} + 10{\zeta}^2)q + ({\zeta}^3 + 108{\zeta}^2 - 513{\zeta}
\\ &+ 808  - \cdots)q^2  + (-64{\zeta}^3 + 808{\zeta}^2 - 2752\zeta +
4016 - \cdots)q^3 + \cdots.
\end{aligned}
$$

The coefficients for $a$ and $b$ can be obtained using Table 1 or
the recursion formulas on page 39 of \cite{ez}. The representation
of $\phi_{2}$ in terms of $a$ and $b$ is:

\begin{center}
$\displaystyle \phi_{2}=\frac{1}{12}a(E_4(\tau)b-E_6(\tau)a)$.
\end{center}

By Theorem 3.2, we have $t^{(2)}(207)=-B(207)$. As $8n-r^2=207$
has a solution $n=29$ and $r=5$, then $B(207)$ is the coefficient
of $q^{29}{\zeta}^5$ which is $-113643$. Thus

\begin{center}
$t^{(2)}(207)=113643 \equiv 0 \bmod 3$.
\end{center}

\end{exa}

\begin{rem} (1) Zagier actually defined $t^{(N)}(d)$ and proved part (ii) of Theorem
3.2 for all $N$ such that $\Gamma_0(N)^{*}$ is of genus 0 (see
Section 8 in \cite{z}). One might be able to prove part (i) of
Theorem 3.2 in the case $N$ is squarefree. If so, then a
congruence, similar to Theorem 1.3, should hold for $t^{(N)}(d)$,
$N$ squarefree. If $N$ is not squarefree, then C. Kim has kindly
pointed out part (i) of Theorem 3.2 does not hold. For example, if
$N=4$ and $d=7$, one can construct $f_{7,4}$ (see the appendix in
\cite{kim2}) and compute that

\begin{center}
$f_{7,4}=q^{-7} -55q + 0q^{4} + 220q^{9} + \cdots.$
\end{center}

\noindent Thus $A(1,7)=-55$. But $B(7)=23$ (see Remark 3.1).

(2) We should note that Theorem 1.3 is an extension of the
simplest case of Theorem 1 in \cite{ao}. Ono and Ahlgren have also
proven congruences for $t(d)$ which involve ramified or inert
primes in quadratic fields. In fact, they prove that a positive
proportion of primes yield congruences for $t(d)$ (see parts (2)
and (3) of Theorem 1 in \cite{ao}). It would be interesting to see
if such congruences hold for $t^{(p)}(d)$ or $t^{(N)}(d)$.

(3) The Monster $\mathbb{M}$ is the largest of the sporadic simple
groups of order

\begin{center}
$2^{46}3^{20}5^{9}7^{6}11^{2}13^{3}17\cdot19\cdot23\cdot29\cdot31\cdot41\cdot47\cdot59\cdot71$
\end{center}

Ogg \cite{ogg} noticed that the primes dividing the order of
$\mathbb{M}$ are exactly those in the set $\frak{S}$. The monster
$\mathbb{M}$ acts on a graded vector algebra $V=V_{-1}
\bigoplus_{n\geq 1} V_n$ (see Frenkel, Lepowsky, and Meurman \cite
{flm} for the construction). For any element $g \in \mathbb{M}$,
let $Tr(g|V_n)$ denote the trace of $g$ acting on $V_n$ for each
$n$. Then $Tr(g|V_{-1})=1$ and $Tr(g|V_n) \in \mathbb Z$ for every
$n \geq 1$. The Thompson series is defined by:

\begin{center}
$T_{g}(z)= q^{-1} + \displaystyle \sum_{n \geq 1} Tr(g|V_n)q^n$.
\end{center}

The authors in \cite{cy} study Thompson series evaluated at
imaginary quadratic arguments, i.e. ``singular moduli'' of
Thompson series. It is possible to define a trace of singular
moduli of Thompson series. A natural question is ``do we have
congruences for these traces?''

\end{rem}

\section*{Acknowledgments}
The author would like to thank Imin Chen, Chang Heon Kim, Ken Ono,
and Noriko Yui for their valuable comments.

\end{document}